\documentclass[graybox]{svmult}
\usepackage{amsrefs}

\usepackage{mathptmx}       % selects Times Roman as basic font
\usepackage{helvet}         % selects Helvetica as sans-serif font
\usepackage{courier}        % selects Courier as typewriter font
\usepackage{makeidx}         % allows index generation
\usepackage{graphicx}        % standard LaTeX graphics tool
\usepackage{multicol}        % used for the two-column index
\usepackage[bottom]{footmisc}% places footnotes at page bottom
\usepackage{amsmath,amssymb,amsfonts}        % AMS math packages
\makeindex             % used for the subject index

\usepackage{tabularx}
\usepackage{makecell}
\usepackage{yfonts}
\newcommand{\R}{\mathbb{R}}

\renewcommand{\epsilon}{\varepsilon}
\renewcommand{\Re}{{\mbox{\textfrak{Re}}}}

\usepackage{url}
\renewcommand{\leq}{\leqslant}
\renewcommand{\le}{\leqslant}

\renewcommand{\ge}{\geqslant}
        
\begin{document}

\title*{Decay estimates in time \\ for classical and anomalous diffusion}

\titlerunning{Classical and anomalous diffusion} 

\author{Elisa Affili, Serena Dipierro and Enrico Valdinoci}

\institute{Elisa Affili \at
Dipartimento di Matematica, Universit\`a degli studi di Milano,
Via Saldini 50, 20133 Milan, Italy, and 
Centre d'Analyse et de Math\'ematique Sociales,
\'Ecole des Hautes \'Etudes en Sciences Sociales,
54 Boulevard Raspail,
75006 Paris, France, \email{elisa.affili@unimi.it}
\and Serena Dipierro \at
Department of Mathematics and Statistics,
University of Western Australia,
35 Stirling Highway,
Crawley WA 6009, Australia, \email{serena.dipierro@uwa.edu.au}
\and Enrico Valdinoci \at
Department of Mathematics and Statistics,
University of Western Australia,
35 Stirling Highway,
Crawley WA 6009, Australia, \email{enrico.valdinoci@uwa.edu.au}
}
\maketitle
\abstract{We present a series of results focused on the decay in time of solutions
of classical and anomalous diffusive equations in a bounded domain.
The size of the solution is measured in a Lebesgue space, and the setting comprises
time-fractional and space-fractional equations and operators of nonlinear type.
We also discuss how fractional operators may affect long-time asymptotics.}

\section{Decay estimates, methods, results and perspectives}

In this note we present some results, recently obtained in~\cite{vespri, elisa},
focused on the long-time behavior of solutions of
evolution equations
which may exhibit anomalous diffusion, caused by either time-fractional
or space-fractional effects (or both). The case of several nonlinear operators
will be also taken into account (and indeed some of the results
that we present are new also for classical diffusion run by nonlinear operators).

The results that we establish give quantitative bounds
on the decay in time of smooth solutions, confined in a smooth bounded
set with Dirichlet data. The size of the solution will be measured
in classical Lebesgue spaces, and we will detect different types of decays according
to the different cases that we take into consideration (the main order
of decay being affected by the
structure of the diffusion in time and by the possible nonlinear character of
the spatial operator).

The evolution equation that we take into account is very general,
and it can be written as an initial datum problem with
homogeneous external Dirichlet condition
of the type
\begin{equation}\label{EQ}
\begin{cases}
\lambda_1\partial_t^\alpha u+\lambda_2\partial_t u+{\mathcal{N}}[u]=0 & {\mbox{ in }}\Omega\times(0,+\infty),\\
u=0 & {\mbox{ in }}(\R^n\setminus\Omega)\times(0,+\infty),\\
u(\cdot,0)=u_0(\cdot) & {\mbox{ in }}\Omega.
\end{cases}\end{equation}
In this setting, $u=u(x,t)$ is a smooth solution of~\eqref{EQ},
$\Omega$ is a bounded set of~$\R^n$ with smooth boundary
(and we are not trying here
to optimize the smoothness assumptions on the solution or on the domain), 
the convex parameters~$\lambda_1$, $\lambda_2\in[0,1]$
are such that~$\lambda_1+\lambda_2=1$, the (possibly nonlinear) operator~${\mathcal{N}}$
acts on the space variable~$x$,
and the time-fractional parameter~$\alpha$ lies
in~$(0,1)$.

Also, in our setting, the symbol~$\partial_t^\alpha$ stands for the so-called
Caputo time-fractional derivative, defined, up to normalizing constants that we omit
for simplicity, by
$$ \partial_t^\alpha v(t):=\frac{d}{dt}\int_0^t \frac{v(\tau)-v(0)}{(t-\tau)^\alpha}\,d\tau.$$
Such a time-fractional derivative naturally arises in many context,
including geophysics~\cite{MR2379269}, neurology~\cite{MR2139665, MR3814763} (see also~\cite{MR3071207} and the references therein) and
viscoelasticity~\cite{MR2351653}, and can be seen as a natural consequence
of classical models of diffusion in highly ramified media such as combs~\cite{comb}.
In addition, from the mathematical point of view,
equations involving the Caputo derivatives can be framed into
the broad line of research devoted to Volterra type integrodifferential operators,
see~\cite{MR2125407, MR3563229}.
\medskip

The operator~${\mathcal{N}}$ in~\eqref{EQ} takes into account the diffusion
in the space variable, and can be either of classical or of fractional type, and concrete
choices will be made in what follows. More precisely, our setting always
comprises, as particular situations, the cases of diffusion driven by the Laplacian or
by the fractional Laplacian, defined by
$$ (-\Delta)^s u(x):=
\int_{\R^n}\frac{u(x)-u(y)}{|x-y|^{n+2s}}\,dy, \qquad \text{with} \ s\in(0,1)$$
where the integral is taken in the principal value sense (to allow cancellations
near the singularity).

In our framework,
we also deal with the case in which~${\mathcal{N}}$
is nonlinear, studying the cases of the classical~$p$-Laplacian
and porous media diffusion (see~\cite{MR1230384,MR2286292})
$$ \Delta_p u^m := {\rm div} (|\nabla u^m|^{p-2}\nabla u^m),
\qquad{\mbox{with $p\in(1,+\infty)$ and $m\in(0,+\infty)$,}}$$
the case of graphical mean curvature, given in formula (13.1) of \cite{giusti}, % !TeX
$$ {\rm div}\left( \frac{\nabla u}{\sqrt{1+|\nabla u|^2}}\right), $$
the case of the fractional $p$-Laplacian (see e.g.~\cite{MR3861716})
\begin{eqnarray*}&&(-\Delta)^s_pu(x):=
\int_{\R^n}\frac{|u(x)-u(y)|^{p-2}(u(x)-u(y))}{|x-y|^{n+sp}}\,dy,\\&&\qquad{\mbox{with $
p\in(1,+\infty)$ and $s\in(0,1)$,}}\end{eqnarray*}
and possibly even the sum of different nonlinear operator of this type, with coefficients $\beta_j>0$, 
$$  \sum_{j=1}^N \beta_j (-\Delta)^{s_j}_{p_j} u, \qquad \text{with} \
		p_j\in(1,+\infty) \ \text{and}  \ s_j\in(0,1),  $$
the case of the anisotropic fractional Laplacian, that is the sum of fractional directional derivatives in the directions of the space $e_j$, given by
$$(-\Delta_{\beta})^{\sigma} u(x)= \sum_{j=1}^{n} \beta_j (-\partial_{x_j}^2)^{\sigma_j} u(x) $$
for $\beta_j>0$, $\beta=(\beta_1, \dots, \beta_n)$ and  $\sigma=(\sigma_1, \dots, \sigma_n)$, where 
$$ (-\partial_{x_j}^2)^{\sigma_j} u(x) = \int_{\R} \frac{u(x)- u(x+\rho e_j)}{\rho^{1+2\sigma_j}} d\rho, $$
considered for example in \cite{anisotropic}. 
The list of possible diffusion operators continues with
two fractional porous media operators (see~\cite{MR2773189, MR2737788})
\begin{eqnarray*}&&
{\mathcal{P}}_{1,s}(u):=(-\Delta)^s u^m
\qquad{\mbox{with $s\in(0,1)$ and $m\in(0,+\infty)$,}}
\\{\mbox{and }}&&{\mathcal{P}}_{2,s}(u):=
-{\rm div}\,(u \nabla {\mathcal{R}}(u)),\qquad{\mbox{ where }}\;{\mathcal{R}}(u)(x):=\int_{\R^n}
\frac{u(y)}{|x-y|^{n-2s}}\,dy \\
&&\qquad\hskip5cm{\mbox{ and }}\;s\in(0,1),
\end{eqnarray*}
the graphical fractional mean curvature operator (see~\cite{MR3331523})
\begin{eqnarray*}&& {\mathcal{H}}^s(u)(x):=\int_{\R^n} F\left(\frac{u(x)-u(x+y)}{|y|}\right)\frac{dy}{|y|^{n+s}},\\&&\qquad\qquad{\mbox{with
$s\in(0,1)$ and }}F(r):=\int_0^r \frac{d\tau}{(1+\tau^2)^{\frac{n+1+s}{2}}},\end{eqnarray*}
the classical Kirchhoff operator for vibrating strings
$$ {\mathcal{K}}(u)(x):=-M\left( \| \nabla u\|^2_{L^2(\Omega)}\right)\;\Delta u(x),$$
and the fractional Kirchhoff operator (see~\cite{MR3120682})
$$ {\mathcal{K}}_s(u)(x):= M\left( \iint_{\R^n\times\R^n}\frac{|u(y)-u(Y)|^2}{|y-Y|^{n+2s}}\,dy\,dY\right)\;
(-\Delta)^s u(x),$$
with~$M:[0,+\infty)\to[0,+\infty)$ nondecreasing and~$s\in(0,1)$. 

The case of complex valued operators
is also considered, in view of a classical (see~\cite{MR0142894})
and fractional (see~\cite{MR3794886}) magnetic
settings, in which we took into account the operators
\begin{eqnarray*}
&& {\mathcal{M}}u:=-(\nabla -i A)^2 u,\\
{\mbox{and }}&&{\mathcal{M}}_s u(x):=\int_{\R^n}\frac{u(x)-e^{i(x-y)A\left(\frac{x+y}{2}\right)}u(y)}{
|x-y|^{n+2s}}\,dy,\qquad{\mbox{with }}\;s\in(0,1),
\end{eqnarray*}
where~$A:\R^n\to\R^n$ represents the magnetic field.\medskip

For further motivations and additional details on these operators, we refer to~\cite{vespri, elisa}:
here we just mention that,
given the general assumptions that we take,
the operator~${\mathcal{N}}$ in~\eqref{EQ} comprises
many cases of interest in both pure and applied mathematics,
with applications in several disciplines,
see for instance~\cite{MR1809268, MR3469920, MR3645874} for detailed discussions
on anomalous diffusion with several applications in different contexts.\medskip

In our setting, we will obtain decay estimates
in suitable Lebesgue spaces~$L^\ell(\Omega)$, for some appropriate exponent~$\ell\ge1$. The typical estimate
that we establish is that all solutions~$u$ of~\eqref{EQ} satisfy
\begin{equation}\label{DEC}
\| u(\cdot,t)\|_{L^\ell(\Omega)}\le C^*\,\Theta(t)\qquad{\mbox{ for all }}t\ge1,
\end{equation}
where~$C^*>0$ depends on the structural assumptions of the problem
(namely on~$\Omega$, $\lambda_1$, $\lambda_2$, $\alpha$,
$\mathcal{N}$, $u_0$ and~$\ell$), and~$\Theta:[1,+\infty)\to(0,+\infty)$
is an appropriate decay function, described here below in concrete situations, possibly depending on another constant $C>0$. The proof
of the decay in~\eqref{DEC} relies on energy estimates, which are in turn
based on suitable Sobolev embeddings that employ the ``parabolic''
structure of the problem, leading to an appropriate ordinary differential
inequality (if~$\lambda_1$ in~\eqref{EQ} is equal to zero), or an appropriate
integral inequality (if~$\lambda_1=1$), or a mixed differential/integrodifferential
inequality (if~$\lambda_1\in(0,1)$),
for the norm-map~$t\mapsto\| u(\cdot,t)\|_{L^\ell(\Omega)}$.% Given appropriate barriers, solutions to suitable differential equations or integral equations, Theorem 1.1 of \cite{vespri} and Theorems 1.1 and 1.2 of \cite{elisa} produce the decay estimates, whose applications are presented in Tables \ref{TA1} and \ref{TA2}. 
The solutions of the equations related to those inequalities are used as barriers and compared to the function $\| u(\cdot,t)\|_{L^\ell(\Omega)}$, as presented in Theorem 1.1 of \cite{vespri} and Theorems 1.1 and 1.2 of \cite{elisa}.

More precisely, the ``elliptic'' character of the spatial diffusive operator
is encoded in an inequality of the type
\begin{equation}\label{MAINI}
\| u(\cdot,t)\|_{L^\ell(\Omega)}^{\ell-1+\gamma}\le C\,\int_\Omega
|u(x,t)|^{\ell-2}\,\Re \Big( \overline{u}(x,t)\,{\mathcal{N}}u(x,t)\Big)\,dx,
\end{equation}
where~$\gamma$ and~$C$ are positive structural constants, ``$\Re$''
denotes the real part and~$\overline{u}$ is the complex conjugate
of~$u$ (in case the problem is set in the reals, the inequality in~\eqref{MAINI}
obviously simplifies). 

We observe that~\eqref{MAINI} becomes more transparent when~$\ell=2$
and~${\mathcal{N}}u=-\Delta u$, with~$u$ real valued: in such a case, after
an integration by parts which takes into account the Dirichlet datum of~$u$, the inequality in~\eqref{MAINI} boils down to the classical Sobolev-Poincar\'e inequality with~$\gamma=1$.

Once the inequality in~\eqref{MAINI} is established for the operator~$\mathcal{N}$
under consideration, one obtains a bound in terms of an ordinary
differential equation, or more generally of a nonlinear integral
equation on the variable~$t$: depending on~$\gamma$ and on the type of time-derivative, this provides
an estimate on the decay in~$t$ of the norm-map~$t\mapsto\| u(\cdot,t)\|_{L^\ell(\Omega)}$, which can be either polynomial or exponential
(in particular, different operators~$\mathcal{N}$ can lead to different
values of~$\gamma$ and therefore to different asymptotics in time
for the solution~$u$).
\medskip

This strategy, suitably adapted to the different situations, applies
to many operators:
the concrete cases that we comprise are listed explicitly in the Tables~\ref{TA1}
and~\ref{TA2}, which
present
the main results achieved in~\cite{vespri, elisa}. For the first table, the theorems cited in the last column are the ones proving \eqref{MAINI} and a decay estimate in the case $\lambda_1=1$ and $\lambda_2=0$ for the operators in their row. Then, combining these results with Theorem 1.1 and 1.2 of \cite{elisa}, the declared estimates trivially follow. 
% ***AGGIUNTA*** 
However, in Table \ref{TA1} for the first time we apply the estimates for the case $\lambda_2=1$ of \cite{elisa} to the operators analyzed in \cite{vespri}, stating the expected decays in a quantitative way.

{\footnotesize
\begin{center}\begin{table}
\begin{tabular}{l | c | c | c | c | c |}
$\,$ & {\bf Operator ${\mathcal{N}}$} & {\bf Values of~$\lambda_1$, $\lambda_2$} & {\bf Range of~$\ell$} & {\bf Decay rate $\Theta$} & {\bf Reference}
\\ \hline \\
\begin{minipage}[l]{0.15\columnwidth}
{\bf Nonlinear classical diffusion 
 }\end{minipage}
& $\Delta_p u^m$ &
$\lambda_1\in(0,1]$, $\lambda_2\in[0,1)$ & $\ell\in[1,+\infty)$ &
$\Theta(t)=\frac1{t^{\frac\alpha{m(p-1)}}}$ & Thm 1.2 \cite{vespri}
\\ \hline \\
\begin{minipage}[l]{0.15\columnwidth}
	{\bf Nonlinear classical diffusion 
		 }\end{minipage}
& $\Delta_p u^m$ {with}  $(m,p)\neq(1,2)$ &
$\lambda_1=0$, $\lambda_2=1$ & $\ell\in[1,+\infty)$ &
$\Theta(t)=\frac1{t^{\frac1{m(p-1)-1}}}$ & Thm 1.2 \cite{vespri}
\\ \hline \\
\begin{minipage}[l]{0.15\columnwidth}
	{\bf Bi-Laplacian }\end{minipage}
& $\Delta_2 u$  &
$\lambda_1=0$, $\lambda_2=1$ & $\ell\in[1,+\infty)$ &
$\Theta(t)=e^{-\frac{t}{C}}$ & Thm 1.2 \cite{vespri}
\\ \hline \\
\begin{minipage}[l]{0.15\columnwidth}
{\bf Graphical mean curvature}\end{minipage}
& $ {\rm div}\left( \frac{\nabla u}{\sqrt{1+|\nabla u|^2}}\right)$ &
$\lambda_1\in(0,1]$, $\lambda_2\in[0,1)$ & $\ell\in[1,+\infty)$ &
$\Theta(t)=\frac1{t^\alpha}$ & Thm 1.5 \cite{vespri}
\\ \hline \\
\begin{minipage}[l]{0.15\columnwidth}
	{\bf Graphical mean curvature}\end{minipage}
& $ {\rm div}\left( \frac{\nabla u}{\sqrt{1+|\nabla u|^2}}\right)$ &
$\lambda_1=0$, $\lambda_2=1$ & $\ell\in[1,+\infty)$ &
$\Theta(t)=e^{-\frac{t}{C}}$ & Thm 1.5 \cite{vespri}
\\ \hline \\
\begin{minipage}[l]{0.15\columnwidth}
{\bf Fractional $p{\mbox{-}}$Laplacian}\end{minipage}
& $(-\Delta)^s_pu$ &
$\lambda_1\in(0,1]$, $\lambda_2\in[0,1)$ & $\ell\in[1,+\infty)$ &
$\Theta(t)=\frac1{t^{\frac\alpha{p-1}}}$ & Thm 1.6 \cite{vespri}
\\ \hline \\
\begin{minipage}[l]{0.15\columnwidth}
	{\bf Fractional $p{\mbox{-}}$Laplacian}\end{minipage}
& $(-\Delta)^s_pu$, $p> 2$ &
$\lambda_1=0$, $\lambda_2=1$ & $\ell\in[1,+\infty)$ &
$\Theta(t)=\frac1{t^{\frac1{p-2}}}$ & Thm 1.6 \cite{vespri}
\\ \hline \\
\begin{minipage}[l]{0.15\columnwidth}
	{\bf Fractional $p{\mbox{-}}$Laplacian}\end{minipage}
& $(-\Delta)^s_pu$, $p\leq 2$ &
$\lambda_1=0$, $\lambda_2=1$ & $\ell\in[1,+\infty)$ &
$\Theta(t)=e^{-\frac{t}{C}}$ & Thm 1.6 \cite{vespri}
\\ \hline \\
\begin{minipage}[l]{0.15\columnwidth} % !TeX Tab1
{\bf Superposition of fractional $p{\mbox{-}}$Laplacians}\end{minipage}
&{ $ \sum_{j=1}^N \beta_j (-\Delta)^{s_j}_{p_j} u$, $\beta_j>0$}&
$\lambda_1\in(0,1]$, $\lambda_2\in[0,1)$ & $\ell\in[1,+\infty)$ &
$\Theta(t)=\frac1{t^{\frac\alpha{p_{\max}-1}}}$ & Thm 1.7 \cite{vespri}
\\ \hline \\
\begin{minipage}[l]{0.15\columnwidth}
	{\bf Superposition of fractional $p{\mbox{-}}$Laplacians}\end{minipage}
& \makecell{$ \sum_{j=1}^N \beta_j (-\Delta)^{s_j}_{p_j} u$, \\ with $\beta_j>0$ and $p_{\max}>2$} &
$\lambda_1=0$, $\lambda_2=1$ & $\ell\in[1,+\infty)$ &
$\Theta(t)=\frac1{t^{\frac1{p_{\max}-2}}}$ & Thm 1.7 \cite{vespri}
\\ \hline \\
\begin{minipage}[l]{0.15\columnwidth}
	{\bf Superposition of fractional $p{\mbox{-}}$Laplacians}\end{minipage}
& \makecell{$ \sum_{j=1}^N \beta_j (-\Delta)^{s_j}_{p_j} u$, \\ with $\beta_j>0$ and $p_{\max}\leq2$}&
$\lambda_1=0$, $\lambda_2=1$ & $\ell\in[1,+\infty)$ &
$\Theta(t)=e^{-\frac{t}{C}}$ & Thm 1.7 \cite{vespri}
\\ \hline \\
\begin{minipage}[l]{0.15\columnwidth}
{\bf Superposition of anisotropic fractional Laplacians}\end{minipage}
& $ \sum_{j=1}^N \beta_j (-\partial_{x_j}^2)^{s_j} u$, $\beta_j>0$&
$\lambda_1\in(0,1]$, $\lambda_2\in[0,1)$ & $\ell\in[1,+\infty)$ &
$\Theta(t)=\frac1{t^{\alpha}}$ & Thm 1.8 \cite{vespri}
\\ \hline \\
\begin{minipage}[l]{0.15\columnwidth}
	{\bf Superposition of anisotropic fractional Laplacians}\end{minipage}
& $ \sum_{j=1}^N \beta_j (-\partial_{x_j}^2)^{s_j} u$, $\beta_j>0$&
$\lambda_1=0$, $\lambda_2=1$ & $\ell\in[1,+\infty)$ &
$\Theta(t)=e^{-\frac{t}{C}}$ & Thm 1.8 \cite{vespri}
\\ \hline \\
\begin{minipage}[l]{0.15\columnwidth} 
{\bf Fractional porous media~I}\end{minipage}
& $ {\mathcal{P}}_{1,s}(u)$&
$\lambda_1\in(0,1]$, $\lambda_2\in[0,1)$ & $\ell\in[1,+\infty)$ &
$\Theta(t)=\frac1{t^{\frac\alpha{m}}}$ & Thm 1.9 \cite{vespri}
\\ \hline \\
\begin{minipage}[l]{0.15\columnwidth}
	{\bf Fractional porous media~I}\end{minipage}
& $ {\mathcal{P}}_{1,s}(u)$, $m>1$ &
$\lambda_1=0$, $\lambda_2=1$ & $\ell\in[1,+\infty)$ &
$\Theta(t)=\frac1{t^{\frac1{m-1}}}$ & Thm 1.9 \cite{vespri}
\\ \hline \\
\begin{minipage}[l]{0.15\columnwidth}
	{\bf Fractional porous media~I}\end{minipage}
& $ {\mathcal{P}}_{1,s}(u)$, $m\leq1$ &
$\lambda_1=0$, $\lambda_2=1$ & $\ell\in[1,+\infty)$ &
$\Theta(t)=e^{-\frac{t}{C}}$ & Thm 1.9 \cite{vespri}
\\ \hline \\
\begin{minipage}[l]{0.15\columnwidth}
{\bf Fractional graphical mean curvature}\end{minipage}
& $ {\mathcal{H}}^s(u)$&
$\lambda_1\in(0,1]$, $\lambda_2\in[0,1)$ & $\ell\in[1,+\infty)$ &
$\Theta(t)=\frac1{t^{\alpha}}$ & Thm 1.10 \cite{vespri}
\\ \hline \\
\begin{minipage}[l]{0.15\columnwidth}
	{\bf Fractional graphical mean curvature}\end{minipage}
& $ {\mathcal{H}}^s(u)$&
$\lambda_1=0$, $\lambda_2=1$ & $\ell\in[1,+\infty)$ &
$\Theta(t)=e^{-\frac{t}{C}}$ & Thm 1.10 \cite{vespri}
\\ \hline
\end{tabular}\medskip
\caption{Results from~\cite{vespri}.}\label{TA1}
\end{table}\end{center}
}

{\footnotesize 
\begin{center}\begin{table}
\begin{tabular}{l | c | c | c | c | c |}
$\,$ & {\bf Operator ${\mathcal{N}}$} & {\bf Values of~$\lambda_1$, $\lambda_2$} & {\bf Range of~$\ell$} & {\bf Decay rate $\Theta$} & {\bf Reference}
\\ \hline \\
\begin{minipage}[l]{0.15\columnwidth}
{\bf Fractional porous media II }\end{minipage}
& ${\mathcal{P}}_{2,s}$ &
$\lambda_1\in(0,1]$, $\lambda_2\in[0,1)$ & $\ell\in[1,+\infty)$ &
$\Theta(t)=\frac1{t^{\frac\alpha{2}}}$ & Thm 1.3 \cite{elisa}
\\ \hline \\
\begin{minipage}[l]{0.15\columnwidth}
{\bf Fractional porous media II }\end{minipage}
& ${\mathcal{P}}_{2,s}$ &
$\lambda_1=0$, $\lambda_2=1$ & $\ell\in[1,+\infty)$ &
$\Theta(t)=\frac1{t}$ & Thm 1.3 \cite{elisa}
\\ \hline \\
\begin{minipage}[l]{0.15\columnwidth}
{\bf Classical Kirchhoff operator}\end{minipage}
& $ {\mathcal{K}}(u)$ with $M(0)>0$ &
$\lambda_1\in(0,1]$, $\lambda_2\in[0,1)$ & $\ell\in[1,+\infty)$ &
$\Theta(t)=\frac1{t^\alpha}$ & Thm 1.4 \cite{elisa}
\\ \hline \\
\begin{minipage}[l]{0.15\columnwidth}
{\bf Classical Kirchhoff operator}\end{minipage}
& $ \begin{matrix}
{\mathcal{K}}(u) {\mbox{ with }} M(t)=bt, \\
b>0{\mbox{ and }}n\le4\end{matrix} $&
$\lambda_1\in(0,1]$, $\lambda_2\in[0,1)$ & $\ell\in[1,+\infty)$ &
$\Theta(t)=\frac1{t^\frac\alpha3}$ & Thm 1.4 \cite{elisa}
\\ \hline \\
\begin{minipage}[l]{0.15\columnwidth}
{\bf Classical Kirchhoff operator}\end{minipage}
& $ \begin{matrix}
{\mathcal{K}}(u) {\mbox{ with }} M(t)=bt, \\
b>0{\mbox{ and }}n\ge5\end{matrix} $&
$\lambda_1\in(0,1]$, $\lambda_2\in[0,1)$ & $\ell\in\left[1,\frac{2n}{n-4}\right)$ &
$\Theta(t)=\frac1{t^\frac\alpha3}$ & Thm 1.4 \cite{elisa}
\\ \hline \\
\begin{minipage}[l]{0.15\columnwidth}
{\bf Classical Kirchhoff operator}\end{minipage}
& $ {\mathcal{K}}(u)$ with $M(0)>0$ &
$\lambda_1=0$, $\lambda_2=1$ & $\ell\in[1,+\infty)$ &
$\Theta(t)=e^{-\frac{t}C}$ & Thm 1.4 \cite{elisa}
\\ \hline \\
\begin{minipage}[l]{0.15\columnwidth}
{\bf Classical Kirchhoff operator}\end{minipage}
& $\begin{matrix}
{\mathcal{K}}(u) {\mbox{ with }}M(t)=bt,\\
b>0\end{matrix}$&
$\lambda_1=0$, $\lambda_2=1$ & $\ell\in[1,+\infty)$ &
$\Theta(t)=\frac1{\sqrt{t}}$ & Thm 1.4 \cite{elisa}
\\ \hline \\
\begin{minipage}[l]{0.15\columnwidth}
{\bf Fractional Kirchhoff operator}\end{minipage}
& $ {\mathcal{K}}_s(u)$ with $M(0)>0$ &
$\lambda_1\in(0,1]$, $\lambda_2\in[0,1)$ & $\ell\in[1,+\infty)$ &
$\Theta(t)=\frac1{t^\alpha}$ & Thm 1.5 \cite{elisa}
\\ \hline \\
\begin{minipage}[l]{0.15\columnwidth}
{\bf Fractional Kirchhoff operator}\end{minipage}
& $ \begin{matrix}
{\mathcal{K}}(u) {\mbox{ with }} M(t)=bt, \\ b>0
{\mbox{ and }}n\le4s\end{matrix} $&
$\lambda_1\in(0,1]$, $\lambda_2\in[0,1)$ & $\ell\in[1,+\infty)$ &
$\Theta(t)=\frac1{t^\frac\alpha3}$ & Thm 1.5 \cite{elisa}\\ \hline \\
\begin{minipage}[l]{0.15\columnwidth}
{\bf Fractional Kirchhoff operator}\end{minipage}
& $ \begin{matrix}
{\mathcal{K}}(u) {\mbox{ with }} M(t)=bt, \\
b>0{\mbox{ and }}n>4s\end{matrix} $&
$\lambda_1\in(0,1]$, $\lambda_2\in[0,1)$ & $\ell\in\left[1,\frac{2n}{n-4s}\right)$ &
$\Theta(t)=\frac1{t^{\frac\alpha3}}$ & Thm 1.5 \cite{elisa}
\\ \hline \\
\begin{minipage}[l]{0.15\columnwidth}
{\bf Fractional Kirchhoff operator}\end{minipage}
& $ {\mathcal{K}}_s(u)$ with $M(0)>0$ &
$\lambda_1=0$, $\lambda_2=1$ & $\ell\in[1,+\infty)$ &
$\Theta(t)=e^{-\frac{t}C}$ & Thm 1.5 \cite{elisa}\\ \hline \\
\begin{minipage}[l]{0.15\columnwidth}
{\bf Fractional Kirchhoff operator}\end{minipage}
& $\begin{matrix} {\mathcal{K}}_s(u)
{\mbox{ with }}M(t)=bt,\\
b>0\end{matrix}$ &
$\lambda_1=0$, $\lambda_2=1$ & $\ell\in[1,+\infty)$ &
$\Theta(t)=\frac1{\sqrt{t}}$ & Thm 1.5 \cite{elisa}
\\ \hline \\
\begin{minipage}[l]{0.15\columnwidth}
{\bf Classical magnetic operator}\end{minipage}
& $ {\mathcal{M}}(u)$ &
$\lambda_1\in(0,1]$, $\lambda_2\in[0,1)$ & $\ell\in[1,+\infty)$ &
$\Theta(t)=\frac1{t^\alpha}$ & Thm 1.6 \cite{elisa}\\ \hline \\
\begin{minipage}[l]{0.15\columnwidth}
{\bf Classical magnetic operator}\end{minipage}
& $ {\mathcal{M}}(u)$ &
$\lambda_1=0$, $\lambda_2=1$ & $\ell\in[1,+\infty)$ &
$\Theta(t)=e^{-\frac{t}C}$ & Thm 1.6 \cite{elisa}
\\ \hline \\
\begin{minipage}[l]{0.15\columnwidth}
{\bf Fractional magnetic operator}\end{minipage}
& $ {\mathcal{M}}_s(u)$ &
$\lambda_1\in(0,1]$, $\lambda_2\in[0,1)$ & $\ell\in[1,+\infty)$ &
$\Theta(t)=\frac1{t^\alpha}$ & Thm 1.7 \cite{elisa}\\ \hline \\
\begin{minipage}[l]{0.15\columnwidth}
{\bf Fractional magnetic operator}\end{minipage}
& $ {\mathcal{M}}_s(u)$ &
$\lambda_1=0$, $\lambda_2=1$ & $\ell\in[1,+\infty)$ &
$\Theta(t)=e^{-\frac{t}C}$ & Thm 1.7 \cite{elisa}\\ \hline 
\end{tabular}\medskip
\caption{Results from~\cite{elisa}.}\label{TA2}\end{table}\end{center}
}

It would be interesting to detect the optimality of the estimates listed in Tables~\ref{TA1} and~\ref{TA2}, and to investigate
other cases of interest as well.
For related decay estimates, see~\cite{MR2454809, MR3294409, MR3563229}.
As a matter of fact, decay estimates for evolutionary problems
are a classical topic of research that has produced
a very abundant, and extremely interesting, literature.
Without aiming at providing an exhaustive list of all the important contributions
on this topic, we mention that:
\begin{itemize}
\item the classical doubly-nonlinear operator~$\Delta_p u^m$ with~$\lambda_2 =1$
has been addressed in~\cite{MR2268115},
\item the classical $1$-Laplace operator has been dealt with in~\cite{MR2033382, KURD, MAZON},
\item decay estimates for the fractional $p$-Laplacian~$(-\Delta_p)^s$
with~$s\in(0,1)$ and~$p>1$ with~$\lambda_2:=1$ have been first established in
Section~6 of~\cite{BCM} (see also~\cite{BCM0}),
\item the case of the fractional $1$-Laplacian, namely~$(-\Delta_p)^s$
with~$s\in(0,1)$ and~$p:=1$ has been treated in~\cite{MAZON},
\item the porous medium equation for $\lambda_2=1$ has been deeply analyzed in \cite{MR3294409},
\item some interesting estimates for the Kirchoff equation are given in \cite{MR2454809},
\item see also~\cite{MR2529737},
where several decay estimates have been obtained
by using integral inequalities.\end{itemize}

We also remark that the interplay between time derivatives
and fractional diffusion produces interesting decay patterns also
in nonlinear equations, see e.g.~\cite{MR3511786}.
Furthermore,
in general, the fractional aspect of the problem can cause significant
differences, as can be observed also from the time decays
of Tables~\ref{TA1} and~\ref{TA2}. For instance, one may notice that
the decay switch from polynomial to exponential in the Kirchhoff equations
when the time-diffusion changes from fractional to classical, and this
independently on the fact that the space diffusion is of classical or fractional type
(roughly speaking, in this context, it is just the character of the diffusion
in time which detects the time decay, regardless the character of
the diffusion in space).
 
We also point out that classical and fractional operators share several
common properties, but they also exhibit structural differences
at a fundamental level. For instance, 
to exhibit an elementary but very interesting feature
in which long-time behaviors are affected by fractional environments, we recall that
fractional diffusion in space,
as modeled by the fractional Laplacian~$(-\Delta)^s$ with~$s\in(0,1)$,
is related to L\'{e}vy-type and $2s$-stable stochastic processes,
and in such case the long jumps of the underlying random walk
causes significant differences with respect to the classical
Brownian motion. In particular, fractional processes
are typically recurrent only in dimension~$1$ and for values
of~$s$ greater or equal to~$1/2$ (being transient in
dimension~$2$ and higher,
and also  in dimension~$1$ for values
of~$s$ smaller than~$1/2$), and this is an important difference
with respect to the case of Gaussian processes,
which are 
recurrent in dimensions~$1$ and~$2$ (and transient
in dimension~$3$ and higher). See
Example~3.5 in~\cite{MR3529858} and the references therein
for a detailed treatment of recurrence and transiency
for  L\'{e}vy-type processes. In this note, in \S~\ref{Y},
we present a very simple, and somewhat heuristic, discussion
of the recurrence and transiency properties related
to the long jump random walks, based on PDE methods
and completely accessible to a broad audience.
\medskip

For a detailed list of other elementary structural differences
between classical and fractional diffusion see also \S~2.1
in~\cite{abatangelo}.

\section{Recurrence and transiency of long jump random processes}\label{Y}

In this section, we discuss a simple PDE approach
to the recurrence of the long jump random walk related to~$(-\Delta)^s$
in dimension~$1$ and for values
of~$s$ greater or equal to~$1/2$ and to its transiency in
dimension~$2$ and higher,
and also  in dimension~$1$ for values
of~$s$ smaller than~$1/2$. The treatment will comprise
the classical case~$s=1$ as well, showing how the structural
differences between the different regimes naturally arise
from a PDE analysis.\medskip

To this end,
for~$s\in(0,1]$, we denote by~${\mathcal{G}}_s(x,t)$ the solution
of the (possibly fractional) heat equation with initial datum
given by the Dirac's Delta, namely
$$ \begin{cases}
\partial_t {\mathcal{G}}_s +(-\Delta)^s {\mathcal{G}}_s=0
& {\mbox{ in }}\R^n\times(0,+\infty),\\
{\mathcal{G}}_s(x,0)=\delta_0(x).
\end{cases}$$
When~$s=1$, we have that such function reduces to the classical
Gauss kernel for the heat flow, namely
$$ {\mathcal{G}}_1(x,t)=
\frac{1}{(4t)^{\frac{n}2}}  e^{-\frac{|x|^2}{4t}}.$$
In general, when~$s\in(0,1)$, the expression of~${\mathcal{G}}_s$
is less explicit, except when~$s=1/2$; in the latter case, it holds that
$$ {\mathcal{G}}_{1/2}(x,t)=\frac{ct}{ (t^2+|x|^2)^{\frac{n+1}{2}}},$$
where~$c>0$ is a normalizing constant -- the need of which lying
in the general mass conservation law
\begin{equation}\label{MASS}
\int_{\R^n} {\mathcal{G}}_s(x,t)\,dx=1,
\end{equation}
see also formula~(2.29) in~\cite{abatangelo}
and page~1363 in~\cite{MR2588782}.
Furthermore, see again page~1363 in~\cite{MR2588782},
we have that
\begin{equation}\label{POSG}
{\mathcal{G}}_s(x,t)>0\qquad{\mbox{for all }}\,(x,t)\in\R^n\times(0,+\infty),\end{equation}
and it enjoys the natural scaling property
\begin{equation}\label{Gscala}
{\mathcal{G}}_s(x,t)=
\frac1{t^{\frac{n}{2s}}}{\mathcal{G}}_s\left(
\frac{x}{t^{\frac{1}{2s}}},1\right).
\end{equation}
See \cite{MR2283957},
\cite{MR2588782}, formulas~(2.41)--(2.45) in~\cite{abatangelo}
and the references therein for a discussion about the fractional
heat kernel and its differences with the classical case.\medskip

For every~$k\in\{1,2,3,\dots\}$ and~$\rho>0$, we define
\begin{equation}\label{qksr} q_k(s,\rho):=\int_{\R^n\setminus B_\rho} {\mathcal{G}}_s(x,k)\,dx.\end{equation}
We observe that
\begin{equation}\label{PROB} 0\le q_k(s,\rho)\le \int_{\R^n} {\mathcal{G}}_s(x,k)\,dx=1,\end{equation}
thanks to~\eqref{MASS}.
Let also
\begin{equation}\label{dqIA}
\begin{split}&
q(s,\rho):=\prod_{k=1}^{+\infty} q_k(s,\rho)\in[0,1]\\ {\mbox{and }}\quad&
q(s):=\lim_{\rho\to0} q(s,\rho).\end{split}
\end{equation}
In view of~\eqref{PROB}, we can consider~$q(s)$
as related to the probability of the stochastic process associated
with the operator~$(-\Delta)^s$ of ``drifting away
without coming back''.
Namely (see~\cite{SEMA}), we know that
$$ \int_A {\mathcal{G}}_s(x,t)\,dx$$
represents the probability that a particle starting at the origin at time~$0$ and following the stochastic process
producing~$(-\Delta)^s$ ends up in the region~$A\subseteq\R^n$
at time~$t$. In this sense, the quantity~$q_k(s,\rho)$
in~\eqref{qksr} represents the probability that this particle
lies outside~$B_\rho$ at time~$k$.

Roughly speaking, for small~$\rho$, a natural Ansatz is to assume these events to be more or less independent from each other: indeed,
in view of~\eqref{Gscala},
using the substitution~$y:=x/\rho$
we have that
\begin{eqnarray*}&&
q_k(s,\rho)=\frac1{k^{\frac{n}{2s}}}\int_{\R^n\setminus B_\rho} 
{\mathcal{G}}_s\left(
\frac{x}{k^{\frac{1}{2s}}},1\right)\,dx=
\frac{\rho^n}{k^{\frac{n}{2s}}}\int_{\R^n\setminus B_1} 
{\mathcal{G}}_s\left(
\frac{\rho y}{k^{\frac{1}{2s}}},1\right)\,dy\\&&\qquad\qquad=
\frac{1}{(k/\rho^{2s})^{\frac{n}{2s}}}\int_{\R^n\setminus B_1} 
{\mathcal{G}}_s\left(
\frac{y}{(k/\rho^{2s})^{\frac{1}{2s}}},1\right)\,dy=
\int_{\R^n\setminus B_1} 
{\mathcal{G}}_s\left(
y,\frac{k}{\rho^{2s}}\right)\,dy
,\end{eqnarray*}
representing the probability of a particle to lie outside~$B_1$
at time~$k/\rho^{2s}$. In view of this, since the time steps~$
k/\rho^{2s}$ are very separated from each other when~$\rho$
is small, we may think that the quantity~$q(s,\rho)$
in~\eqref{dqIA}
is a good approximation of the probability that the particle
does not lie in~$B_1$ 
in all the time steps~$
k/\rho^{2s}$, as well as 
a good approximation of the probability that the particle
does not lie in~$B_\rho$ 
in all the time steps~$
k\in\{1,2,3,\dots\}$. In this heuristics,
the case in which the quantity~$q(s)$
in~\eqref{dqIA} is equal to~$0$ indicates
that the particle will come back 
infinitely often
to its original position at the origin
in integer times (with probability~$1$); conversely,
the case in which the quantity~$q(s)$
in~\eqref{dqIA} is equal to $1$ indicates
that the particle will return
to its original position at the origin
in integer times
only with probability zero.

In this sense, computing~$q(s)$ gives an interesting indication
of the recurrence properties of the associated stochastic process,
and, in our case, this calculation can be performed as follows.
First of all, we notice that
$$ \inf_{x\in B_1}{\mathcal{G}}_s(x,1):=\iota_s>0,$$
thanks to~\eqref{POSG}, and
$$ \sup_{x\in\R^n}{\mathcal{G}}_s(x,1):=\mu_s<+\infty.$$
As a consequence, recalling~\eqref{MASS} and~\eqref{Gscala},
if~$\rho\in(0,1]$
\begin{eqnarray*}&& p_k(s,\rho):=1-q_k(s,\rho)=
\int_{B_\rho} {\mathcal{G}}_s(x,k)\,dx\\&&\qquad=
\frac1{k^{\frac{n}{2s}}}
\int_{B_\rho} {\mathcal{G}}_s\left(
\frac{x}{k^{\frac{1}{2s}}},1\right)\,dx\in \left[ 
\frac{\iota_s\,|B_\rho|}{ k^{\frac{n}{2s}} },
\frac{\mu_s\,|B_\rho|}{ k^{\frac{n}{2s}}}
\right].
\end{eqnarray*}
This gives that
$$ 
q_k(s,\rho)=1-p_k(s,\rho)\in\left[1-\frac{C\rho^n}{ k^{\frac{n}{2s}}}, 1-\frac{c\rho^n}{ k^{\frac{n}{2s}}}\right],$$
for some~$C>c>0$, depending only on~$n$ and~$s$, and accordingly
\begin{equation}\label{qro}\begin{split}&
\log q(s,\rho)=\log\left(\prod_{k=1}^{+\infty} q_k(s,\rho)\right)\\&\qquad
=\sum_{k=1}^{+\infty}
\log q_k(s,\rho)\in\left[
\sum_{k=1}^{+\infty}\log\left(1-\frac{C\rho^n}{ k^{\frac{n}{2s}}}\right)
,\,
\sum_{k=1}^{+\infty}\log\left(1-\frac{c\rho^n}{ k^{\frac{n}{2s}}}\right)
\right].\end{split}\end{equation}
Also, for a fixed~$C_0\in(0,+\infty)$, the convergence
of the series
$$ \sum_{k=1}^{+\infty}\log\left(1-\frac{C_0\rho^n}{ k^{\frac{n}{2s}}}\right)$$
can be reduced to that of the series
$$ -\sum_{k=1}^{+\infty}\frac{C_0\rho^n}{ k^{\frac{n}{2s}}},$$
and consequently
$$ \sum_{k=1}^{+\infty}\log\left(1-\frac{C_0\rho^n}{ k^{\frac{n}{2s}}}\right)=\begin{cases}
-C_1\rho^n & {\mbox{ if }}n>2s,\\
-\infty& {\mbox{ if }}n\le2s,
\end{cases}$$
for some~$C_1>0$ depending only on~$n$, $s$ and~$C_0$.
This and~\eqref{qro} lead to
\begin{eqnarray*}
&& \log q(s,\rho)=-\infty \qquad{\mbox{ if }}n\leq2s,\\
&& \log q(s,\rho)\in[-C_2\rho^n,-C_3\rho^n]\qquad{\mbox{ if }}n>2s,\end{eqnarray*}
for some~$C_2>C_3>0$ depending only on~$n$ and $s$.

Therefore,
\begin{eqnarray*}
&& q(s,\rho)= 0\qquad{\mbox{ if }}n\le2s,\\
&& q(s,\rho)\in[e^{-C_2\rho^n},e^{-C_3\rho^n}]\qquad{\mbox{ if }}n>2s.\end{eqnarray*}
Taking the limit as~$\rho\to0$, we thereby find that
\begin{equation}\label{DRI}
\begin{split}
& q(s)= 0\qquad{\mbox{ if }}n\le2s,\\
& q(s)=1\qquad{\mbox{ if }}n>2s.\end{split}\end{equation}
When~$s=1$, we can write~\eqref{DRI} as
\begin{equation*}
\begin{split}
& q(s)=0\qquad{\mbox{ if }}n\in\{1,2\},\\
& q(s)= 1\qquad{\mbox{ if }}n\ge3,
\end{split}\end{equation*}
that is, in our framework, the classical random walk ``comes back to the initial'' point in dimensions~$1$ and~$2$, and it ``drifts
away forever'' in dimension~$3$ and higher.

When~$s\in(0,1)$, the situation is different, since~\eqref{DRI}
produces the alternative
\begin{equation*}
\begin{split}
& q(s)=0\qquad{\mbox{ if $n=1$ and $s\in[1/2,1)$}},\\
& q(s)= 1\qquad{\mbox{ if $n\ge2$, and also if~$n=1$ and~$s\in(0,1/2)$.}}
\end{split}\end{equation*}
That is,
in our setting, the fractional random walk ``comes back to the initial'' point only in dimensions~$1$ and
only if the fractional parameter is above a certain threshold
(namely~$s\ge1/2$). Conversely, the
fractional random walk
``drifts
away for ever'' already in dimension~$2$,
and even in dimension~$1$ if the fractional parameter is too
small (namely~$s<1/2$). 

See e.g.~\cite{MR3738798}
and the
references therein for a comprehensive treatment of recurrence
and transiency of general stochastic processes.

\begin{acknowledgement}
This work has been supported by
the Australian Research Council Discovery Project 170104880 NEW ``Nonlocal
Equations at Work''. Part of this work was carried out
on the occasion of a very pleasant visit of the first author
to the University of Western Australia, which we thank
for the warm hospitality.
The authors are members of INdAM/GNAMPA.\end{acknowledgement}

\end{document}